%% file: brezinski_arxiv_2.tex
\documentclass[9pt,a4paper]{article}                     
\usepackage{algorithm2e}
\usepackage[latin1]{inputenc}
\usepackage[T1]{fontenc}
\usepackage{bezier}
\usepackage{enumerate}
\usepackage{dsfont}
\usepackage{mathptmx}      
\usepackage{amsfonts}
\usepackage{amsmath}
\usepackage{amssymb}
\usepackage{amsbsy}
\usepackage{amscd}
\usepackage{amsopn}
\usepackage{amssymb}
\usepackage{amstext}
\usepackage{amsxtra}
\usepackage{graphicx,psfrag}
\usepackage{latexsym}
\usepackage{multirow}
\newtheorem{thm}{Theorem}
\newtheorem{lem}{Lemma}

\def\squareforqed{\hbox{\rlap{$\sqcap$}$\sqcup$}}
\def\qed{\ifmmode\else\unskip\quad\fi\squareforqed}
\def\smartqed{\def\qed{\ifmmode\squareforqed\else{\unskip\nobreak\hfil
\penalty50\hskip1em\null\nobreak\hfil\squareforqed
\parfillskip=0pt\finalhyphendemerits=0\endgraf}\fi}}
\smartqed

\DeclareMathAlphabet{\mathcal}{OMS}{cmsy}{m}{n}

\newcommand{\wabs}[1]{\left|#1\right|}

\newcommand{\wc}{\mathds{C}}

\newcommand{\wcal}[1]{\mathcal{#1}}

\newcommand{\wdfdx}[3]{\frac{\partial {#1}}{\partial {#2}}\!\wlr{{#3}}}

\newcommand{\wfc}[2]{{#1}\!\left(#2\right)}

\newcommand{\wi}[1]{\wrm{i}}

\newcommand{\wlr}[1]{\left( #1 \right)}

\newcommand{\wn}{\mathds N}

\newcommand{\wnorm}[1]{\left\| #1 \right\|}

\newcommand{\wpade}{Pad\'{e}}

\newcommand{\wref}[1]{(\ref{#1})}

\newcommand{\wrm}[1]{\mathrm{#1}}

\newcommand{\wseq}[2]{{\left\{ {#1}_{#2}, \ {#2} \in \wn \right\}}}

\newcommand{\wset}[1]{{\left\{ #1 \right\}}}

\newcommand{\wvec}[1]{\mathbf{#1}}

\newcommand{\wz}{\mathds Z}

\begin{document}

\title{The divergence of the barycentric \wpade{} interpolants
\thanks{This work was supported by grant \#2013/10916-2, S\~{a}o Paulo Research Foundation (FAPESP).}
}
\author{Walter F. Mascarenhas
\thanks{Instituto de Matem\'{a}tica e Estat\'{i}stica, Universidade de S\~{a}o Paulo,
\           Cidade Universit\'{a}ria, Rua do Mat\~{a}o 1010, S\~{a}o Paulo SP, Brazil. CEP 05508-090
              Tel.: +55-11-3091 5411, Fax: +55-11-3091 6134,  walter.mascarenhas@gmail.com}
}
\maketitle

\begin{abstract}
We explain that, like the usual \wpade{} approximants,
the barycentric \wpade{} approximants proposed recently by Brezinski and Redivo-Zaglia
can diverge. More precisely, we show that for every polynomial
$\wfc{P}{z}$ there exists a function $\wfc{g}{z} = \sum_{n = 0}^\infty c_n z^n$, with $c_n$ arbitrarily small,
such that the sequence of barycentric \wpade{} approximants of
$\wfc{f}{z} = \wfc{P}{z} + \wfc{g}{z}$ do not converge uniformly
in any subset of $\wc{}$ with a non-empty interior.
\end{abstract}


\input intro_2
\input main_proof_2
\input easy_lemmas_2
\input proof_of_lemmas_2

\input references_2

\end{document}

%% file: intro_2.tex
\section{Introduction}
In the recent article \cite{BREZINSKI}, Claude Brezinski and Michela Redivo-Zaglia proposed
a barycentric version of \wpade{} approximation and illustrated its effectiveness in practice.
In the conclusion of their article they asked whether their approximants converge in theory.
In this article we explain that, like the usual \wpade{} approximants, there are entire functions
for which the barycenctric \wpade{} approximants do not converge uniformly in any
subset of $\wc{}$ with a non-empty interior.

In the barycentric approach to \wpade{} approximation
proposed by Brezinski and Redivo--Zaglia, given $n + 1$ distinct interpolation
points $x_{n,m} \in \wc{} - \wset{0}$, we define
\begin{equation}
\label{def_pq}
\wfc{p_n}{z} := \sum_{m = 0}^n \frac{w_{n,m} \wfc{f}{x_{n,m}}}{z - x_{n,m}} \hspace{1cm} \wrm{and} \hspace{1cm}
\wfc{q_n}{z} := \sum_{m = 0}^n \frac{w_{n,m}}{z - x_{n,m}},
\end{equation}
with weights $w_{n,m}$ chosen so that $\wfc{f}{z} \wfc{q_n}{z} = \wfc{p_n}{z} + \wfc{O}{z^n}$.
The resulting barycentric approximants $\wfc{p_n}{z}/\wfc{q_n}{z}$ interpolate $\wfc{f}{z}$ at the
points $x_{n,m}$ and match its first $n-1$ derivatives
at $z = 0$. Of course, usual \wpade{} approximants with the same degrees of freedom would
match more derivatives at $z = 0$ and
the barycentric approach exchanges these derivatives by the interpolation at the points $x_{n,m}$.

Given a polynomial $\wfc{P}{z}$,
interpolation nodes $X = \wset{x_{n,m}, \, n \in \wn{}, \,  0 \leq m \leq n} \subset \wc{} - \wset{0}$,
with $x_{n,m} \neq x_{n,k}$ for $m \neq k$, and a set
$\wseq{\alpha}{k} \subset \wc{} - X$, we explain how to build functions of the form
\begin{equation}
\label{main_f}
\wfc{f}{z} = \wfc{P}{z} + \sum_{m = 0}^\infty c_m z^m,
\end{equation}
with $c_m$ arbitrarily small, and  indexes $\wseq{n}{k}$,
such that $f$'s barycentric \wpade{} approximant of degree $n_k$ has a pole
arbitrarily close to $\alpha_k$. This
shows that the poles of $f$'s barycentric approximants can form a dense subset of $\wc{}$.
In this case, the sequence of approximants do not converge uniformly to $\wfc{f}{z}$ in any set with
a non-empty interior.

In formal terms, we prove the following theorem. The $\alpha_k$ represent
complex numbers to which the poles of the approximants will be arbitrarily close. We
show that $\alpha_k$ can be arbitrarily chosen, as long as they do not coincide
with the interpolation nodes. Once the $\alpha_k$ are chosen, the theorem guarantees the
existence of approximants with poles very close to them.
\begin{thm}
\label{thm_main}
Consider
\begin{enumerate}[(1)]
\item A set $X = \wset{x_{n,m}, \, n \in \wn{}, \,  0 \leq m \leq n} \subset \wc{} - \wset{0}$
 with $x_{n,m} \neq x_{n,k}$ for $m \neq k$.
\item A sequence $\wseq{\epsilon}{m}$ of small positive numbers.
\item A sequence $\wseq{n}{k}$ of indexes with $n_{k + 1} > 2 n_k$.
\item A sequence $\wseq{\alpha}{k} \subset \wc{} - X$.
\end{enumerate}
For every polynomial $\wfc{P}{z}$ with $\wfc{\wrm{degree \ of}}{P} < n_0$
there exists a set $\wseq{\pi}{k} \subset \wc{} - X$ and coefficients $\wseq{c}{m}$ such that
the function $\wfc{f}{z}$ in \wref{main_f} is entire and
\begin{enumerate}[(i)]
\item $c_m = 0$ for $m < n_0$ and $\wabs{c_m} \leq \epsilon_m$ for $m \geq n_0$.
\item $\wabs{\pi_k - \alpha_k} \leq \epsilon_{n_k}$ for $k \in \wn{}$.
\item For all $k \in \wn{}$, $\pi_k$ is a pole of the barycentric \wpade{} approximant
of $\wfc{f}{z}$ with nodes $x_{n_k,m}$.
\end{enumerate}
\end{thm}

In  section \ref{sec_main_proof} we prove Theorem \ref{thm_main}. Our proof uses lemmas which
are stated in section \ref{sec_lemmas} and proved in section \ref{sec_proof_lemmas}. We suggest
that, at first, the reader follows the proof of the general theorem accepting
the lemmas as true results. Unfortunately, the proof is very technical. In order to
motivate our technical arguments,
in the next section we present them in an informal way and sketch
an algorithm to compute the coefficients $c_m$ and the poles $\pi_k$ mentioned
in the main theorem. We hope that by reading this informal section
the reader will grasp the overall structure of the proof and will no be
distracted by the unavoidable technical details.

\section{An informal description of the proof of Theorem \ref{thm_main}}
\label{section_informal}
We can think of the proof of Theorem \ref{thm_main} as an algorithm to
compute the coefficients $c_m$ in such way that we guarantee the existence
of the poles $\pi_k$ near the points $\alpha_k$. Unfortunately, this
algorithm would not work in practice because it considers a sequence of
approximants for which the number of interpolation nodes grows exponentially.
It also uses infinite sequences, of which most
terms would underflow in finite precision arithmetic.
Moreover, rounding errors would change things completely. In other
words, our examples are neither robust nor practical.

The fact that our algorithm to construct the function $\wfc{f}{z}$
would not work in practice is an evidence of a deeper problem with this
article and the theoretical analysis of the convergence of
numerical algorithms in general. As we have shown in other instances
\cite{MEA,MEB,MEC,MED,MEE}, several algorithms which work well in practice
are vulnerable to theoretical examples like the ones presented here.
Therefore, the algorithm by Brezinski and Redivo--Zaglia may be quite
adequate for all the degrees of the approximants one usually considers
in practice (but this would be the subject of another article.)

It is our opinion that our examples tell more about the inadequacy of
the asymptotic analysis of the algorithms than they tell about the
inadequacy of the algorithms themselves. In fact, articles like this
one only show how far we are from an adequate theory to explain
the behavior of these algorithms for large (but finite) problem instances
in finite precision arithmetic.

The difference between theory and practice affects not only the convergence
of numerical algorithms; it affects applied mathematics in general. In
Statistics the discussion of this topic is quite old,
as one can notice by reading the article \cite{MALLOWS} by C. Mallows,
the commentaries after it and its references. In Computer Science the discussion is illustrated
in the first 15 minutes of the provocative talk by Alan Kay at OOPSALA97 \cite{KAY}
and the last minutes of Donald Knuth's talk at Google \cite{KNUTH}.
The difference of opinions of outstanding scholars like Peter Huber and Brad Efron illustrated
in Mallows's article and the opposing views of Kay and Knuth
show that there are no easy answers regarding
the interplay of theory and practice in applied mathematics
(in a broad sense.)

That said, we now present an informal algorithm to compute the
coefficients $c_m$ in Theorem \ref{thm_main}, under the assumption
that we are using exact arithmetic and can estimate constants
for which our usual theories provide only existence proofs
(like the radius of convergence for Newton's method for a
$C^1$ function for which we do not know how to bound the derivative
or its inverse.)
To keep things simple, we assume that we are concerned with
the polynomial $\wfc{P}{z} = 0$ and take $n_k = 3^{k}$.

The first step is to write the entire function $f$ in theorem \ref{thm_main} as
\begin{equation}
\label{intuitive_f}
\wfc{f}{z} = \sum_{k = 0}^\infty \mu_k z^{3^k} \sum_{m = 0}^{3^k} d_{k,m} z^m,
\end{equation}
where the $\mu_k$ are free parameters to be determined by our algorithm
and the $d_{k,m}$ are carefully chosen constants which depend
on the interpolation nodes and the $\alpha_k$, but which
do not depend on the $\mu_k$. The algebraic expressions defining the
constants $d_{k,m}$ are relevant for the technical details but, once the reader believes
in our claims about them, they do not matter much for the overall understanding
of the proof.

We build the $\mu_k$ one by one, by induction. However, we must be careful because the
location of the poles of the barycentric Padé approximant with nodes
\[
X_k := \wset{x_{3^k,m}, m = 0,\dots, 3^k}
\]
of the function $f$ in \wref{intuitive_f} will be also influenced by the
$\mu_n$'s for $n > k$, which we do not know at the $k$th step.
The idea then is to bound the $\mu_n$ for $n > k$ so that they do not influence
much the location of the poles of the previous approximants.
Therefore, instead of building only
the sequence $\mu_k$ we also build a family $\beta_{k,j}$ of bounds such that
if, for $n > k$,  $\wabs{\mu_n} \leq \beta_{k,n}$ then the location of the
poles determined at step $k$ will not be significantly affected by the
$\mu_n$ with $n > k$. In order to achieve this goal the $\mu_k$ and the
$\beta_{k,n}$ must decay very rapidly with $k$, so rapidly in fact
that they would underflow in finite precision arithmetic.

In idealized terms, the algorithm can be summarized as follows:
\begin{enumerate}
\item Start with $\beta_{-1,m} = 1$ for $m \in \wn{}$.
\item For $j = 0, 1, \dots, $ do
\item Choose $\mu_j$ with $\wabs{\mu_j} \leq \beta_{j-1,j}$ such that
the barycentric Padé approximant of
\[
\wfc{f_j}{z} = \sum_{k = 0}^j \mu_k z^{3^k} \sum_{m = 0}^{3^k} d_{k,m} z^m
\]
with nodes $X_j = \wset{x_{3^j,m}, \ m = 0,\dots, 3^j}$ has a pole near $\alpha_j$
(the constants $d_{k,m}$ are chosen so that this is possible.)
\item Use the complex version of the implicit function theorem to define
$\beta_{j,n} \leq \beta_{j-1,n}$ for $n \in \wn{}$
such that if $\wabs{\mu_n} \leq \beta_{j,n}$ for $n > k$ then the
the barycentric Padé approximant of the function
$f$ in \wref{intuitive_f} also has a pole near $\alpha_j$.
\item goto 2.
\end{enumerate}

In principle, one could try to apply the same procedure to
analyze other versions of Padé approximants. However, this
may not work because we may not be able to adapt step 3.
We were able to build our examples because we found
constants $d_{m,k}$ with an important property,
which may be specific to the barycentric Padé approximants:
the $d_{m,k}$ are such that the algebraic formulae to compute
the barycentric Padé approximants imply the existence
of the poles near the $\alpha_j$'s. In some sense,
for these $d_{m,k}$, the method causes its own demise,
because the same equations that ensure
its degree of approximation at the origin and the
interpolation at the remaining nodes lead to
the poles near the $\alpha_j$.

Unfortunately, there are many technical details involved in turning the
informal arguments above into a theorem. We tried to find a simpler
way than our proof of Theorem \ref{thm_main} and our lemmas to achieve
this goal, but we failed.

%% file: main_proof_2.tex
\section{Proof of Theorem \ref{thm_main}}
\label{sec_main_proof}
Let us start by defining the terms we use.
We are concerned with sub-sequences with indexes $n_k$
of the sequence of barycentric \wpade{}
approximants. For $n = n_k$, we interpolate at distinct points $x_{n_k,0}$, $x_{n_k,1}$, $\dots$, $x_{n_k,n_k}$
and define
\[
t_{k,m} := x_{n_k,m}.
\]
The Vandermonde matrix $\wvec{V}_k$ corresponding to Lagrange interpolation at
$t_{k,m}$ is
\begin{equation}
\label{def_vk}
\wvec{V}_k =
\left(
\begin{array}{ccccc}
1      & t_{k,0}   & t_{k,0}^2 & \dots  & t_{k,0}^{n_k} \\
1      & t_{k,1}   & t_{k,1}^2 & \dots  & t_{k,1}^{n_k} \\
1      & t_{k,2}   & t_{k,2}^2 &        & \vdots \\
\vdots & \vdots    &           & \ddots & \vdots \\
1      & t_{k,n_k} &  \dots    & \dots  & t_{k,n_k}^{n_k}
\end{array}
\right),
\hspace{1.0cm} \wrm{with} \hspace{1.0cm}
v_{k,i,j} := t_{k,i}^j.
\end{equation}
(Our matrices have indexes $(i,j)$, with $0 \leq i, j \leq n_k$, and we denote
$\wlr{\wvec{V}_k}_{i,j}$ by $v_{k,i,j}$.)
The weights for usual barycentric interpolation at $t_{k,m}$ are
\begin{equation}
\label{def_lambda}
\lambda_{k,m} := \frac{1}{\prod_{i \neq m} \wlr{t_{k,i} - t_{k,m}}}.
\end{equation}
We also use the vectors $\wvec{a}_k \in \wc{}^{1 + n_k}$ with entries
\begin{equation}
\label{def_ak}
a_{k,m} := \frac{1}{\alpha_k - t_{k,m}},
\end{equation}
and $\wvec{B}_k$  is the $\wlr{1 + n_k} \times \wlr{1 + n_k}$
 diagonal matrix which has $a_{k,i}$ in its diagonal:
\begin{equation}
\label{def_bk}
b_{k,i,i} = a_{k,i} \hspace{1cm} \wrm{and} \hspace{1cm} b_{k,i,j} = 0 \ \ \wrm{for} \ \ i \neq j.
\end{equation}
The $\wlr{1 + n_k}$-dimensional vector $\wvec{e}_k$ has entries
\begin{equation}
\label{def_ek}
e_{k,i} := 0 \ \ \wrm{for} \ \ 0 \leq i < n_k \hspace{1cm} \wrm{and} \hspace{1cm} e_{k,n_k} := 1.
\end{equation}
The coefficients $c_m$ of the function $\wfc{f}{z}$ in \wref{main_f} are defined in terms of the
vectors
\begin{equation}
\label{def_dk}
\wvec{d}_k := \wvec{V}^{-1}_k \wvec{a}_k,
\end{equation}
and a sequence $\wseq{\sigma}{k}$:
\begin{eqnarray}
\label{def_cm_p}
c_m & := & 0 \ \ \hspace{1.22cm} \wrm{for} \ \ m < n_0, \\
\label{def_cm_0}
c_m & := & 0 \ \
\hspace{1.22cm} \wrm{for} \ k \in \wn{} \hspace{0.2cm}
\wrm{and} \hspace{0.2cm}  2 n_k < m < n_{k+1}, \\
\label{def_cm}
c_{m} & := & \sigma_k d_{k, m - n_k} \hspace{0.2cm} \wrm{for} \ k \in \wn{} \hspace{0.2cm}
\wrm{and}
\hspace{0.2cm} n_k \leq m \leq 2 n_k,
\end{eqnarray}
so that
\begin{equation}
\label{def_fg}
\wfc{f}{z} = \wfc{f}{z;\sigma} := \wfc{P}{z} + \sum_{k = 0}^\infty \sigma_k z^{n_k} \sum_{m = 0}^{n_k} d_{k,m} z^m.
\end{equation}

Let us define
\begin{equation}
\label{def_rk}
r_k := 1 + \max_{0 \leq m \leq k, \ \ 0 \leq j \leq n_m} \wabs{t_{m,j}},
\end{equation}

\[
\tau_k := \frac{\min_{\ 0 \leq m \leq 2 n_{k + 1}} \epsilon_m}{\wlr{1 +
\sum_{m = 0}^{2 n_{k+1}} \epsilon_m} r_{k+1}^{2 n_{k+1}} \wlr{1 + \wnorm{\wvec{d}_{k+1}}_1} \wlr{1 + n_{k+1}}!}.
\]
Note that $0 < \tau_k < 1$ and,  if $\sigma_k$ in \wref{def_cm} is such that $0 < \sigma_k \leq \tau_{k-1}$, then
$\wabs{c_m} < \epsilon_m$ for $n_k \leq m \leq 2 n_k$.
Moreover,
\begin{equation}
\label{def_s}
\chi := \sum_{k = 0}^\infty r_k^{2 n_k} \wnorm{\wvec{d}_k}_1 \tau_k < \infty,
\end{equation}
and the series in \wref{def_fg} converges for all $z$ when $0 \leq \sigma_k \leq \tau_{k-1}$ for all $k$.

Finally, by perturbing $\alpha_k$, we can assume that
\begin{equation}
\label{good_alpha}
\sum_{m = 0}^{n_k} \frac{\lambda_{k,m}}{\alpha_k - t_{k,m}} \neq 0
\hspace{1cm} \wrm{and} \hspace{1cm}
\sum_{m = 0}^{n_k} \frac{ \lambda_{k,m} \, t_{k,m}^{n_k}}{\alpha_k - t_{k,m}} \neq 0.
\end{equation}

We are now ready to prove Theorem \ref{thm_main}.

{\bf Proof of Theorem \ref{thm_main}.} We prove the following:
\begin{quote}
Main claim: There exists $\wseq{\mu}{k}$ so that the
function $\wfc{f}{z} = \wfc{f}{z;\mu}$ in \wref{def_fg}
satisfies the requirements of Theorem \ref{thm_main}.
\end{quote}

In order to verity the main claim, we build
$\wset{\mu_k, \ k \in \wz{}}$
and $\wset{\rho_k, \ k \in \wz{}}$ such that:
\begin{enumerate}[(a)]
\item $\rho_{m} := \mu_m := 1$  for $m < 0$.
\item $0 < \rho_{m+1} \leq \rho_m$ for all $m \in \wz{}$.
\item $0 < \mu_m \leq \rho_{m-1} \tau_{m-1}$  for $m = 0,1,2,\dots$.
\item
Let $m \in \wn{}$ and $\wseq{\sigma}{h}$ be such that
\footnote{The numbers $\rho_{m} \tau_{h-1}$ are the bound $\beta_{k,j}$ mentioned in the
 argument in section \ref{section_informal}.}
\begin{enumerate}[(i)]
\item $\sigma_h = \mu_h$ for $0 \leq h \leq m$,
\item $0 < \sigma_h \leq \rho_m \tau_{h-1}$ for $h > m$,
\end{enumerate}
then there exists $\wfc{\xi}{\sigma} \in \wc{}$ such that $\wabs{\wfc{\xi}{\sigma} - \alpha_m} \leq \tau_m$
and the barycentric \wpade{} approximant $\wfc{p_{n_m}}{z}/\wfc{q_{n_m}}{z}$ of the function
$\wfc{f}{z} = \wfc{f}{z; \sigma}$ in \wref{def_fg} satisfies
\begin{equation}
\label{hard_part}
\wfc{p_{n_m}}{\wfc{\xi}{\sigma}} \neq 0 \hspace{0.7cm} \wrm{and} \hspace{0.7cm} \wfc{q_{n_m}}{\wfc{\xi}{\sigma}} = 0.
\end{equation}
\end{enumerate}
 
The existence of $\rho_k$ and $\mu_k$ satisfying (a)--(d) verifies the main claim
because, for each $m \in \wn{}$, we can apply item (d) to
$\sigma = \mu$ and conclude that there exists $\pi_m = \wfc{\xi}{\sigma}$ as required by Theorem \ref{thm_main}.

We have already defined $\mu_k$  and $\rho_k$ for $k < 0$ and
the items (a)--(d) above hold for negative $m = k < 0$.
We now assume that $k\geq 0$ and we have defined $\mu_m$ and $\rho_{m}$ for $m < k$ and the items
(b)--(d) hold for such $m$, and
define $\mu_k$ and $\rho_k$ such that (b)--(d) holds for $m < k + 1$.
By the induction principle, this defines $\mu_k$ and $\rho_k$ for all $k \in \wz{}$.

For $\mu_0, \mu_1, \dots, \mu_k$, consider the function
\begin{equation}
\label{def_fk}
\wfc{f_k}{z} := \wfc{f_k}{z;\mu} :=
\wfc{P}{z} + \sum_{h = 0}^{k - 1} \mu_h z^{n_h} \sum_{m = 0}^{n_h} d_{h,m} z^m +
\mu_k z^{n_k} \sum_{m = 0}^{n_k} d_{k,m} z^m.
\end{equation}
Let $\wseq{\sigma}{m}$ be such that
\begin{equation}
\label{def_sigma}
\sigma_h = \mu_h \ \wrm{for} \ \ 0 \leq h \leq k \hspace{0.5cm} \wrm{and} \hspace{0.5cm}
0 < \sigma_h \leq \rho_k \tau_h \ \ \wrm{for} \ \ h > k.
\end{equation}
The Lemmas \ref{lem_decomp}, \ref{lem_non_deg} and  \ref{lem_gk} show that
the barycentric \wpade{} approximant for the function $\wfc{f}{z} = \wfc{f}{z; \sigma}$
in \wref{def_fg} for $\sigma$ in \wref{def_sigma} is defined by
matrices $\wvec{Y}$, $\wvec{U}$ and $\wfc{\wvec{S}}{\sigma}$ and
weights $\wfc{\wvec{w}}{\sigma} \neq 0$  with
\begin{equation}
\label{def_wsig}
\wlr{\wvec{Y} + \mu_k \wvec{U} + \wfc{\wvec{S}}{\sigma}} \wfc{\wvec{w}}{\sigma} = 0,
\end{equation}
and, for $\chi$ in \wref{def_s},
\begin{equation}
\label{small_s_of_sigma}
\wnorm{\wfc{\wvec{S}}{\sigma}}_2 < \rho_k \wlr{1 + n_k}\chi.
\end{equation}
These lemmas show that there exist $\mu_k \in (0, \rho_{k-1} \tau_{k-1})$ and $\wvec{v} \in \wc{}^{1 + n_k}$ such that:
\begin{equation}
\label{vm_non_null}
v_m \neq 0 \ \ \wrm{for} \ \ 0 \leq m \leq n_k,
\end{equation}
\begin{equation}
\wfc{\wrm{rank \ of}}{\wvec{Y} + \mu_k \wvec{U}} = n_k \hspace{1cm} \wrm{and} \hspace{1cm} \wlr{\wvec{Y} + \mu_k \wvec{U}} \wvec{v} = 0,
\end{equation}
\begin{equation}
\label{main_good_pnz}
\sum_{m = 0}^{n_k} \frac{v_m}{\wlr{\alpha_k - t_{k,m}}^2} \neq 0 \hspace{1cm} \wrm{and} \hspace{1cm}
\sum_{m = 0}^{n_k} \frac{v_m \wfc{f_k}{t_{k,m}}}{\alpha_k - t_{k,m}} \neq 0.
\end{equation}
(To verify the second inequality in \wref{main_good_pnz}, take $\kappa = -\sum_{m = 0}^{n_k} v_m \wfc{f_{k-1}}{t_{k,m}}/\wlr{\alpha_k - t_{k,m}}$
in Lemma \ref{lem_gk}.)
Since $\wvec{Y} + \mu_k \wvec{U}$ is a $n_k \times \wlr{1 + n_k}$ matrix with rank $n_k$,
there exists $\zeta_0 \in \wlr{0,\tau_k}$  such that
\begin{equation}
\label{rank_cond}
\wnorm{\delta \wvec{M}}_2 \leq \zeta_0 \Rightarrow \wvec{Y} + \mu_k \wvec{U} + \delta \wvec{M} \ \ \wrm{has \ rank \ } n_k.
\end{equation}
By continuity and \wref{main_good_pnz}, there exists $\zeta_1 \in \wlr{0, \min \wset{\zeta_0, \wabs{v_0}, \wabs{v_1}, \dots, \wabs{v_{n_k}}}}$
such that
\begin{equation}
\label{perturb_pnz}
\max \wset{ \wnorm{\delta \wvec{v}}_2, \wabs{\delta \alpha}, \wabs{\delta \!f_m}} \leq \zeta_1 \Rightarrow
\sum_{m = 0}^{n_k} \frac{\wlr{v_m + \delta \! v_m}
\wlr{\wfc{f_k}{t_{k,m}} + \delta \! \! f_m}}
{\alpha_k + \delta \alpha - t_{k,m}} \neq 0.
\end{equation}

Lemma \ref{lem_decomp} shows that the entries in the first row of $\wvec{Y}$
are all zero and that
the first equation in the system $\wlr{\wvec{Y} + \mu_k \wvec{U}} \wvec{v} = 0$ can be written as
$\wfc{\eta}{\alpha_k,\wvec{v}} = 0$, for
\[
\wfc{\eta}{\alpha,\wvec{v}} := \sum_{m = 0}^{n_k} \frac{v_m}{\alpha - t_{k,m}} = 0.
\]
Equation \wref{main_good_pnz} shows that the function $\wfc{\eta}{\alpha,\wvec{v}}$ has partial derivative
\[
\wdfdx{\!\eta}{\alpha}{\alpha_k,\wvec{v}} = \sum_{m = 0}^{n_k} \frac{v_m}{\wlr{\alpha_k - t_{k,m}}^2} \neq 0.
\]
Since $\wfc{\eta}{\alpha_k,\wvec{v}} = 0$, the (complex) implicit function theorem shows that
there exists $\zeta_2 \in (0,\zeta_1)$ such that if
$\wnorm{\delta \! \wvec{v}}_2 < \zeta_2$ then there exists
$\wfc{\theta}{\delta \! \wvec{v}} \in \wc{}$ with
\begin{equation}
\label{def_theta}
\wabs{\wfc{\theta}{\delta \! \wvec{v}}} \leq \zeta_1 \hspace{0.4cm} \wrm{and} \hspace{0.4cm}
\wfc{\eta}{\alpha + \wfc{\theta}{\delta \! \wvec{v}},
\wvec{v} + \delta \! \wvec{v}} = \sum_{m = 0}^{n_k} \frac{v_m + \delta \! v_m}{\alpha_k + \wfc{\theta}{\delta \! \wvec{v}} - t_{k,m}} = 0.
\end{equation}
Since $\wvec{Y} + \mu_k \wvec{U}$ has rank $n_k$, there exists $\zeta_3 \in (0,\zeta_2)$ such that
if $\wnorm{\delta \wvec{M}}_2 \leq \zeta_3$ then there exists
$\wfc{\kappa}{\delta \wvec{M}}$ with
\begin{equation}
\label{def_kappa_delta}
\wfc{\kappa}{\delta \wvec{M}} \leq \zeta_2 \hspace{1cm} \wrm{and} \hspace{1cm}
\wlr{\wvec{Y} + \mu_k \wvec{U} + \delta \wvec{M}} \wlr{ \wvec{v} + \wfc{\kappa}{\delta \wvec{M}}} = 0.
\end{equation}

We claim that by considering $\chi$ in \wref{def_s} and taking
\begin{equation}
\label{def_rhok}
\rho_k := \min \wset{\rho_{k-1}, \frac{\zeta_3}{\wlr{1 + n_k} \chi}}
\end{equation}
and the $\mu_k$ above we satisfy the requirement (d) on $\rho_k$ and $\mu_k$
for $m = k$, and we end this proof validating this claim.
In fact, let $\wseq{\sigma}{h}$ be a sequence satisfying \wref{def_sigma}.
Equations  \wref{small_s_of_sigma} and \wref{def_rhok}  show that $\wnorm{\wfc{\wvec{S}}{\sigma}}_2 \leq \zeta_3$
and \wref{rank_cond} implies that the matrix $\wvec{Y} + \mu_k \wvec{U} + \wfc{\wvec{S}}{\sigma}$
has rank $n_k$. Therefore, the space of solutions $\wfc{\wvec{w}}{\sigma}$ of
\wref{def_wsig} has
dimension one. Equation \wref{def_kappa_delta} shows that
\[
\tilde{\wvec{w}} := \wvec{v} + \wfc{\kappa}{\wfc{\wvec{S}}{\sigma}}
\]
is a solution of \wref{def_wsig}. It follows that all solutions $\wfc{\wvec{w}}{\sigma}$
of \wref{def_wsig} are of the form $\gamma \tilde{\wvec{w}}$, with $\gamma \in \wc{}$.
Since all these solutions lead to
the same approximant ($\gamma$ cancels out), the approximants
are defined by $\tilde{\wvec{w}}$.

Equation \wref{def_kappa_delta} shows that $\delta \wvec{v} = \wfc{\kappa}{\wfc{\wvec{S}}{\sigma}}$
is such that $\wnorm{\delta \wvec{v}} \leq \zeta_2$ and leads to $\wfc{\theta}{\delta \wvec{v}}$
satisfying \wref{def_theta}. Since
\wref{def_theta} is equivalent to $\wfc{q_{n_k}}{\wfc{\xi}{\sigma}} = 0$ for
$\wfc{\xi}{\sigma} := \alpha_k + \wfc{\theta}{\delta \wvec{v}}$, we have
verified the last condition in \wref{hard_part}. Moreover,
$\wabs{\wfc{\xi}{\sigma} - \alpha_k} = \wabs{\wfc{\theta}{\delta \wvec{v}}} < \zeta_1 < \tau_k$.

Consider $z$ with $\wabs{z} < r_k$, with $r_k$ in \wref{def_rk} and $f_k$ in \wref{def_fk}.
Since $\wabs{\sigma_h} \leq \rho_k \tau_h$ for $h > k$, equations \wref{def_s} and \wref{def_rhok}
show that $\wfc{\delta\!f}{z} := \wfc{f}{z;\sigma} - \wfc{f_k}{z}$ satisfies
\[
\wabs{\wfc{\delta\!f}{z}} = \wabs{\wfc{f}{z;\sigma} - \wfc{f_k}{z}}
\leq \sum_{h = k+1}^\infty \wabs{\sigma_h}  r_h^{n_h} \sum_{m = 0}^{n_h} \wabs{d_{h,m}} r_h^m \leq
\]
\[
\leq \rho_k \sum_{h = k+1}^\infty \tau_h  r_h^{2 n_h} \wnorm{\wvec{d}_{h}}_1 \leq
\rho_k \wlr{1 + n_k} \chi \leq \zeta_3.
\]
Therefore $\wabs{\delta \! f_m} = \wabs{\wfc{\delta \! f_k}{t_{m,k}}} \leq \zeta_3$  for $0 \leq m \leq n_k$.
Since, for $\wfc{f}{z;\sigma}$ in \wref{def_fg},
\[
\wfc{f}{t_{k,m}} = \wfc{f_k}{t_{k,m};\sigma} + \wfc{\delta \! f_k}{t_{m,k}},
\]
equation  \wref{perturb_pnz} shows that
$\wfc{p_{n_k}}{\wfc{\xi}{\sigma}} = \wfc{p_{n_k}}{\alpha_k + \wfc{\theta}{\delta \wvec{v}}} \neq 0$.
Therefore, we have verified the first condition in \wref{hard_part} and we are done.
\qed{}.

%% file: easy_lemmas_2.tex
\section{Lemmas}
\label{sec_lemmas}

\begin{lem}
\label{lem_system}
For $R > 0$, suppose that $\sum_{m = 0}^\infty \wabs{c_m} R^m < \infty$  and consider distinct
points $x_{n,0},  \dots, x_{n,n}$ with $0 < \wabs{x_{n,m}} < R$.
The functions $\wfc{p_n}{z}$ and $\wfc{q_n}{z}$ in \wref{def_pq}  yield the $n$-th
degree barycentric \wpade{} approximant for $\wfc{f}{z} = \sum_{m = 0}^\infty c_m z^m$
if and only if, for $0 \leq i < n$,
\begin{equation}
\label{eq_system}
\sum_{j = 0}^n \wlr{\sum_{k = n - i}^\infty c_k x_{n,j}^{k - n + i}} w_{n,j} = 0.
\end{equation}
\end{lem}

\begin{lem}
\label{lem_decomp} For the coefficients $c_m$ in \wref{def_cm_p}--\wref{def_cm},
there are matrices $\wvec{Y}$, $\wvec{U}$ and $\wfc{\wvec{S}}{\sigma}$
with dimension $n_k \times \wlr{1 + n_k}$ such that
$\wvec{w} \in \wc{}^{1 + n_k}$ satisfies \wref{eq_system} if and only if
\[
\wlr{\wvec{Y} + \sigma_k \wvec{U} + \wfc{\wvec{S}}{\sigma}} \wvec{w} = 0,
\]
and
\begin{enumerate}[(1)]
\item $\wvec{Y}$ does not depend on $\sigma_m$ for $m \geq k$.
\item All the entries in the first row of $\wvec{Y}$ are equal to zero, i.e.,
$y_{0,j} = 0$.
\item $\wvec{U}$ has entries
\begin{equation}
\label{def_uij}
u_{i,j} = \frac{t_{k,j}^{i}}{\alpha_k - t_{k,j}},
\hspace{0.7cm} \wrm{and, \ in \ particular } \ \
u_{0,j} = \frac{1}{\alpha_k - t_{k,j}}.
\end{equation}
\item If $\sigma_m \leq \epsilon \tau_m$ for all $m > k$ and $\wseq{\tau}{m}$ satisfies
\wref{def_s} then
\[
\wnorm{\wfc{\wvec{S}}{\sigma}}_2 \leq \epsilon \wlr{1 + n_k} \chi.
\]
\end{enumerate}
\end{lem}

\begin{lem}
\label{lem_non_deg}
Let $\wvec{U}$ be as in Lemma \ref{lem_decomp}.
For every matrix $\wvec{M}$ with dimension $n_k \times \wlr{1 + n_k}$,
there exists a finite set $\wcal{E}$ such that if
$\epsilon \not \in \wcal{E}$ then the matrix
$\wvec{M} + \epsilon \wvec{U}$ has rank $n_k$ and there
exists a vector $\wfc{\wvec{v}}{\epsilon}$ such that
$\wlr{\wvec{M} + \epsilon \wvec{U}} \wfc{\wvec{v}}{\epsilon} = 0$, and
for $0 \leq m \leq n_k$,
\begin{enumerate}[(1)]
\item  $\wfc{v_m}{\epsilon} \neq 0$ and
$\wfc{v_m}{\epsilon}$ is a rational function of $\epsilon$.
\item For the vector $\Lambda_k$ with entries $\lambda_{k,m}$ in \wref{def_lambda}
and $\wvec{B}_k$ in \wref{def_bk},
\begin{equation}
\label{lim_v}
\lim_{\epsilon \rightarrow \infty} \wfc{\wvec{v}}{\epsilon} = \wvec{B_k}^{-1} \Lambda_{k}.
\end{equation}
\end{enumerate}
\end{lem}

\begin{lem}
\label{lem_gk}
For $\wvec{d}_k$ in \wref{def_dk}, consider a constant $\kappa \in \wc{}$ and the polynomial
\begin{equation}
\label{def_gk}
\wfc{g_k}{z} := z^{n_k} \sum_{j = 0}^{n_k} d_{k,j} z^j,
\end{equation}
If $\alpha_k$ satisfies \wref{good_alpha}
and $\wfc{\wvec{v}}{\epsilon}$ is a vector whose coordinates are rational
functions of $\epsilon$ and satisfy \wref{lim_v}, then there exists a finite
set $\wcal{E}$ such that if $\epsilon \not \in \wcal{E}$ then
\begin{equation}
\label{good_pnk}
\sum_{m = 0}^{n_k} \frac{
\wfc{v_m}{\epsilon}}
{\wlr{\alpha_k - t_{k,m}}^2} \neq 0
\hspace{1cm} \wrm{and} \hspace{1cm}
\epsilon \sum_{m = 0}^{n_k} \frac{
\wfc{v_m}{\epsilon} \wfc{g_k}{ t_{k,m} }}
{\alpha_k -  t_{k,m}} \neq \kappa.
\end{equation}
\end{lem} 

%% file: proof_of_lemmas_2.tex
\section{Proofs of the lemmas}
\label{sec_proof_lemmas}

{\bf Proof of Lemma \ref{lem_system}.} If $\wabs{z} < \min_{\ 0 \leq j \leq n} \wabs{x_{n,j}}$ then
equation \wref{def_pq} yields
\[
\wfc{p_n}{z} = \sum_{j = 0}^n \frac{w_{n,j}}{z - x_{n,j}} \sum_{k = 0}^\infty c_k x_{n,j}^k
 = \sum_{k = 0}^\infty c_k \wlr{\sum_{j = 0}^n \frac{w_{n,j} x_{n,j}^k}{z - x_{n,j}}} =
 \]
 \[
= \sum_{k = 0}^\infty c_k \wlr{\sum_{j = 0}^n \frac{w_{n,j} x_{n,j}^k}{x_{n,j}}
\frac{1}{\frac{z}{x_{n,j}}  - 1}} =
\]
\begin{equation}
\label{p_exp}
 - \sum_{k = 0}^\infty c_k \wlr{\sum_{j = 0}^n \frac{w_{n,j} x_{n,j}^k}{x_{n,j}}
 \sum_{h = 0}^\infty z^h x_{n,j}^{-h}} =
- \sum_{h = 0}^\infty \wlr{\sum_{k = 0}^\infty c_k \wlr{\sum_{j = 0}^n w_{n,j} x_{n,j}^{k-h-1}}} z^h.
\end{equation}
Moreover,
\[
\wfc{q_n}{z} = \sum_{j=0}^n \frac{w_{n,j}}{z - x_{n,j}} =
\sum_{j=0}^n \frac{w_{n,j}}{x_{n,j}} \frac{1}{\frac{z}{x_{n,j}} - 1}
= - \sum_{j=0}^n \frac{w_{n,j}}{x_{n,j}} \sum_{h = 0}^\infty \frac{z^h}{x_{n,j}^h} =
\]
\begin{equation}
\label{q_exp}
= - \sum_{h = 0}^\infty \wlr{\sum_{j=0}^n w_{n,j} x_{n,j}^{-\wlr{h+1}}} z^h .
\end{equation}
Equation \wref{q_exp} shows that
\begin{equation}
\label{fq}
\wfc{f}{z} \wfc{q_n}{z} = - \sum_{h = 0}^{n-1}
\wlr{
\sum_{k = 0}^h c_k
\wlr{\sum_{j=0}^n w_{n,j} x_{n,j}^{k - h - 1}}
}
z^h + \wfc{O}{z^n}.
\end{equation}

Combining \wref{p_exp} with \wref{fq} we get that
$\wfc{f}{z} \wfc{q_n}{z} = \wfc{p_n}{z} + \wfc{O}{z^n}$
if and only if, for $0 \leq h < n$,
\[
\sum_{k = 0}^\infty c_k \wlr{\sum_{j = 0}^n w_{n,j} x_{n,j}^{k-h-1}} =
\sum_{k = 0}^h c_k \wlr{\sum_{j=0}^n w_{n,j} x_{n,j}^{k - h - 1}}.
\]
Subtracting the right-hand-side from the left-hand-side in this equation we obtain
\[
\sum_{j = 0}^n \wlr{ \sum_{k = h + 1}^\infty  c_k x_{n,j}^{k-h-1}} w_{n,j} = 0,
\]
and replacing $h$ by $n - i - 1$ in the equation above we obtain \wref{eq_system}.
\qed{}

{\bf Proof of Lemma \ref{lem_decomp}.} Equations \wref{def_vk} and \wref{def_dk} show that
\begin{equation}
\label{akj}
\sum_{m = 0}^{n_k} d_{k,m} t_{k,j}^{m} = \wlr{\wvec{V}_k \wvec{d}_k}_j = a_{k,j}.
\end{equation}
Given $0 \leq i < n_k$, we write $\wfc{P}{z} = \sum_{h = 0}^{n-1} p_h z^h$ and define
\begin{eqnarray}
\nonumber
A_i & := & \wset{  0   \leq h < n_0} \cap \wset{h \geq n_k - i}, \\
\nonumber
B_i & := & \wset{n_0   \leq h < n_k} \cap \wset{h \geq n_k - i},
\end{eqnarray}
$\gamma_h := p_h$ for $h \in A_i$ and $\gamma_h := c_h$ for $h \in B_i$.
Equations \wref{def_cm_p}--\wref{def_cm}, \wref{def_uij} and \wref{akj}
show that, for $0 \leq i < n_k$, we have
\[
\sum_{h = n_k - i}^\infty c_h x_{n_k,j}^{h - n_k + i} =
y_{i,j} + \tilde{u}_{i,j} + \wfc{s_{i,j}}{\sigma},
\]
with
\begin{eqnarray}
\label{def_yij}
y_{i,j} & := & \sum_{h \in A_i \cup B_i} \gamma_h t_{k,j}^{h - n_k + i},  \\
\nonumber
\tilde{u}_{i,j} & := & \sum_{h = n_k}^{2 n_k}  c_h t_{k,j}^{h - n_k + i} =
\sigma_k t_{k,j}^i \sum_{m = 0}^{n_k}
d_{k,m} t_{k,j}^{m} = \sigma_k t_{k,j}^i a_{k,i} = \sigma_k u_{i,j},\\
\nonumber
\wfc{s_{i,j}}{\sigma} & := & \sum_{h = 2 n_k + 1}^{\infty}  c_h x_{n_k,j}^{h - n_k + i}
= \sum_{l = k + 1}^\infty \sigma_l \sum_{m = 0}^{n_l}  d_{l,m} t_{k,j}^{n_l - n_k + i + m}.
\end{eqnarray}
Therefore, the system of equations \wref{eq_system} can be written as
$\wlr{\wvec{Y} + \sigma_k \wvec{U} + \wfc{\wvec{S}}{\sigma}} \wvec{w} = 0$, for the
matrices $\wvec{Y}$ and $\wfc{\wvec{S}}{\sigma}$ with entries $y_{i,j}$ and $\wfc{s_{i,j}}{\sigma}$
above and $u_{i,j}$ in \wref{def_uij}.

Note that $y_{i,j}$ does not depend on $\sigma_m$ for $m \leq k$.
When $i = 0$ we have $A_i \cup B_i = \emptyset$
and, as a result, $y_{0,j} = 0$. Thus, the $y_{ij}$ in \wref{def_yij} satisfy
items (1) and (2) in Lemma \ref{lem_decomp}. Moreover, if $0 \leq \sigma_m \leq \epsilon \tau_m$ then,
for $0 \leq i < n_k$, \wref{def_rk} and \wref{def_s} imply that
\[
\wabs{\wfc{s_{ij}}{\sigma}} \leq
\sum_{l = k + 1}^\infty \sigma_l \sum_{m = 0}^{n_l}  \wabs{d_{l,m}}
\wabs{t_{k,j}}^{n_l - n_k + i + m}
\leq \epsilon  \sum_{l = k + 1}^\infty \tau_l \wnorm{\wvec{d}_l}_1 r_{k}^{2 n_l} \leq \epsilon \chi.
\]
Therefore,
\[
\wnorm{\wfc{\wvec{S}}{\sigma}}_2 \leq
\sqrt{\sum_{0 \leq i,j \leq n_k} \wabs{\wfc{s_{ij}}{\sigma}}^2} \leq \epsilon \wlr{1 + n_k} \chi
\]
and we are done.
\qed{}

{\bf Proof of Lemma \ref{lem_non_deg}.} Let $\tilde{\wvec{M}}$ be the matrix obtained
by adding a null $n_k$-th row to $\wvec{M}$ and let $\tilde{\wvec{U}}$ the matrix we obtain
by adding to $\wvec{U}$ the $n_k$-th row with entries
\[
u_{n_k,j} = \frac{t_{k,j}^{n_k}}{\alpha_k - t_{k,j}}.
\]
The matrix $\tilde{\wvec{U}}$ can be recast as
\begin{equation}
\label{u_tilde}
\tilde{\wvec{U}} = \wvec{V}_k^t \, \wvec{B}_k,
\end{equation}
for $\wvec{V}_k$ in \wref{def_vk} and $\wvec{B}_k$ in \wref{def_bk}.
Thus, $\tilde{\wvec{U}}$ is non-singular and the determinant of
the matrix $\wfc{\wvec{N}}{\epsilon} := \tilde{\wvec{M}} + \epsilon \tilde{\wvec{U}}$
is a polynomial $\wfc{Q}{\epsilon}$. This polynomial is not identically zero,
because the non-singularity of $\tilde{\wvec{U}}$ implies that
$\lim_{\epsilon \rightarrow \infty} \wabs{\wfc{Q}{\epsilon}} = +\infty$.
Therefore, there exists only a finite set of $\epsilon$s for which
$\wfc{\wvec{N}}{\epsilon}$ is singular.
We define $\wcal{E}_{-1}$ as the union of this finite set with $\wset{0}$.

Given $\epsilon \not \in \wcal{E}_1$ and $\wvec{e}_k$ in \wref{def_ek}, the vector
\[
\wfc{\wvec{v}}{\epsilon} := \epsilon \wfc{\wvec{N}}{\epsilon}^{-1} \wvec{e}_k
\]
satisfies $\wlr{\wvec{M} + \epsilon \wvec{U}} \wfc{\wvec{v}}{\epsilon} = 0$
and its coordinates are rational functions of $\epsilon$.
Moreover,
\[
\wfc{\wvec{v}}{\epsilon} = \epsilon \wlr{\tilde{\wvec{M}} + \epsilon \tilde{\wvec{U}}}^{-1} \wvec{e}_k
= \wlr{\frac{1}{\epsilon} \tilde{\wvec{M}} + \tilde{\wvec{U}}}^{-1} \wvec{e}_k.
\]
Therefore, \wref{u_tilde} yields
\[
\lim_{\epsilon \rightarrow \infty}  \wfc{\wvec{v}}{\epsilon} =  \tilde{\wvec{U}}^{-1} \wvec{e}_k
= \tilde{\wvec{v}} :=   \wvec{B}^{-1}_k \, \wvec{V}^{-t}_k \, \wvec{e}_k.
\]
Cramer's rule, Laplace's expansion and equation \wref{def_ek} show that
\begin{equation}
\label{eq_detv}
\tilde{v}_m = \wlr{-1}^{n_k + m} \wlr{\alpha_k - t_{k,j}}
\frac{\wfc{\det}{\wvec{W}_{k,m}}}{\wfc{\det}{\wvec{V}_k}},
\end{equation}
where $\wvec{W}_{k,m}$ is the matrix obtained by the removal of 
the $m$-th row and last column of $\wvec{V}_k$. 
$\wvec{V}_k$ and $\wvec{W}_{k,m}$ are Vandermonde matrices and,  therefore,
\[
\wfc{\det}{\wvec{V}_k} = \prod_{0 \leq i < j \leq n} \wlr{t_{k,j} - t_{k,i}}
\hspace{0.5cm} \wrm{and} \hspace{0.5cm}
\wfc{\det}{\wvec{W}_{k,m}} =  \prod_{0 \leq i < j \leq n, \ i,j \neq m} \wlr{t_{k,j} - t_{k,i}}.
\]
The equation above and equations \wref{def_lambda} and \wref{eq_detv} imply that
\[
\tilde{v}_m = \frac{\wlr{\alpha_k - t_{k,j}} \wlr{-1}^{n_k + m}}{\wlr{ \prod_{0 \leq j < m}
\wlr{t_{k,j} - t_{k,m}}}  \wlr{ \prod_{m < i \leq n_k} \wlr{t_{k,m} - t_{k,i}}}} =
\]
\[
= \frac{\alpha_k - t_{k,j}}{\prod_{j \neq m}\wlr{t_{k,j} - t_{k,m}}}
= \wlr{\alpha_k - t_{k,j}} \lambda_{k,j},
\]
and we have verified \wref{lim_v}.

Finally, for every $0 \leq i \leq n_k$ , $\wfc{v_i}{\epsilon}$ is
a rational function of $\epsilon$ and the last paragraph shows that this rational function
does not vanish for large $\epsilon$. This implies that there exists a finite set
$\wcal{E}_i$ such that if $\epsilon \not \in \wcal{E}_i$ then $\wfc{v_i}{\epsilon} \neq 0$.
We complete this proof by taking $\wcal{E} := \bigcup_{i = -1}^{n_k} \wcal{E}_i$.
\qed{}

{\bf Proof of Lemma \ref{lem_gk}.}
Let us show that there exist a finite set $\wcal{E}_1$ such that
if $\epsilon \not \in \wcal{E}_1$ then the first inequality in
\wref{good_pnk} holds. Equations \wref{good_alpha} and
\wref{lim_v} show that the rational function of $\epsilon$ given by
\[
\wfc{\mu}{\epsilon} :=
\sum_{m = 0}^{n_k} \frac{\wfc{v_m}{\epsilon}}
{\wlr{\alpha_k - t_{k,m}}^2}
\]
satisfies
\[
\lim_{\epsilon \rightarrow \infty} \wfc{\mu}{\epsilon}
= \sum_{m = 0}^{n_k} \frac{ \lambda_{k,m}}{\alpha_k -  t_{k,m}} \neq 0.
\]
This implies that the finite set $\wcal{E}_1$ mentioned above exists.

We now prove that there exist a finite set $\wcal{E}_2$ such that
if $\epsilon \not \in \wcal{E}_2$ then the second inequality in
\wref{good_pnk} holds. The definitions of $\wvec{V}_k$, $\wvec{a}_k$ and $\wvec{d}_k$
in \wref{def_vk}, \wref{def_ak}, \wref{def_dk} and \wref{def_gk} yield
\[
\label{eq_gk}
\wfc{g_k}{ t_{k,m}} =  t_{k,m}^{n_k} \sum_{j = 0}^{n_k} d_{k,j}  t_{k,m}^j =
t_{k,m}^{n_k} \sum_{j = 0}^{n_k}  t_{k,m}^j d_{k,j} =
\frac{ t_{k,m}^{n_k}}{\alpha_k -  t_{k,m}}.
\]
This implies that the vector $\wvec{h}$ with coordinates $h_m = \wfc{g_k}{ t_{k,m}}$ satisfies
\begin{equation}
\label{vec_g}
\wvec{h} = \wvec{B}_k \wvec{V}_k \wvec{e}_k,
\end{equation}
for $\wvec{B}_k$ in \wref{def_bk} and $\wvec{e}_k$ in \wref{def_ek}.
Consider the function
\[
\wfc{\gamma}{\epsilon} := \frac{1}{\epsilon} \wlr{
\epsilon \sum_{m = 0}^{n_k} \frac{\wfc{v_m}{\epsilon} \wfc{h_k}{ t_{k,m}}}
{\alpha_k -  t_{k,m}} - \kappa} =
\frac{1}{\epsilon} \wlr{\epsilon \wfc{\wvec{v}}{\epsilon}^t \wvec{B}_k \wvec{h} - \kappa}.
\]
Equations \wref{good_alpha}, \wref{lim_v} and \wref{vec_g} imply that
\[
\lim_{\epsilon \rightarrow \infty} \wfc{\gamma}{\epsilon} =
\Lambda_k^t \wvec{B}_k \wvec{V}_k \wvec{e_k}
= \sum_{m = 0}^{n_k} \frac{ \lambda_{k,m} \, t_{k,m}^{n_k}}{\alpha_k - t_{k,m}} \neq 0,
\]
and Lemma \ref{lem_gk} follows from the observation that $\wfc{\gamma}{\epsilon}$ is
rational function of $\epsilon$.
\qed{}